\documentclass[12pt]{article}
\usepackage{graphicx}
\usepackage{amsmath,amsthm,amssymb,enumerate}
\usepackage{euscript,mathrsfs}
\usepackage[left=2cm,right=2cm,top=3.5cm,bottom=3.5cm]{geometry}
\usepackage{color}

\usepackage{soul}

\catcode`\@=11 \@addtoreset{equation}{section}

\catcode`\@=12

\allowdisplaybreaks

\newtheorem{Theorem}{Theorem}[section]
\newtheorem{Proposition}[Theorem]{Proposition}
\newtheorem{Lemma}[Theorem]{Lemma}
\newtheorem{Corollary}[Theorem]{Corollary}

\theoremstyle{definition}
\newtheorem{Definition}[Theorem]{Definition}

\newtheorem{Remark}[Theorem]{Remark}

\newcommand{\bTheorem}[1]{
\begin{Theorem} \label{T#1} }
\newcommand{\eT}{\end{Theorem}}

\newcommand{\bProposition}[1]{
\begin{Proposition} \label{P#1}}
\newcommand{\eP}{\end{Proposition}}

\newcommand{\bLemma}[1]{
\begin{Lemma} \label{L#1} }
\newcommand{\eL}{\end{Lemma}}

\newcommand{\bCorollary}[1]{
\begin{Corollary} \label{C#1} }
\newcommand{\eC}{\end{Corollary}}

\newcommand{\bRemark}[1]{
\begin{Remark} \label{R#1} }
\newcommand{\eR}{\end{Remark}}

\newcommand{\bDefinition}[1]{
\begin{Definition} \label{D#1} }
\newcommand{\eD}{\end{Definition}}

\newcommand{\Del}{\Delta_x}

\newcommand{\vud}{\vu_\delta}

\newcommand{\bfphi}{\boldsymbol{\varphi}}

\newcommand{\vuB}{\vu_B}
\newcommand{\vrB}{\vr_B}

\newcommand{\bFormula}[1]{
\begin{equation} \label{#1}}
\newcommand{\eF}{\end{equation}}

\newcommand{\Ov}[1]{\overline{#1}}

\newcommand{\DC}{C^\infty_c}

\newcommand{\aleq}{\stackrel{<}{\sim}}
\newcommand{\Ds}{\mathbb{D}_x}

\newcommand{\vr}{\varrho}

\newcommand{\vu}{\vc{u}}

\newcommand{\vrd}{\vr_\delta}

\newcommand{\vc}[1]{{\bf #1}}

\newcommand{\Div}{{\rm div}_x}
\newcommand{\Grad}{\nabla_x}
\newcommand{\Dt}{\frac{\rm d}{{\rm d}t}}

\newcommand{\dx}{\,{\rm d} {x}}

\newcommand{\dt}{\,{\rm d} t }

\newcommand{\dxdt}{\dx  \dt}
\newcommand{\intO}[1]{\int_{\Omega} #1 \ \dx}

\newcommand{\intT}{\int_{\mathcal{T}^1}}

\newcommand{\intTO}[1]{\int_{\mathcal{T}^1} \int_{\Omega} #1 \ \dxdt}

\newcommand{\D}{{\rm d}}
\newcommand{\ep}{\varepsilon}

\def\softd{{\leavevmode\setbox1=\hbox{d}%
          \hbox to 1.05\wd1{d\kern-0.4ex{\char039}\hss}}}
\definecolor{Cgrey}{rgb}{0.85,0.85,0.85}
\definecolor{Cblue}{rgb}{0.50,0.85,0.85}
\definecolor{Cred}{rgb}{1,0,0}
\definecolor{fancy}{rgb}{0.10,0.85,0.10}

\newcommand\Cbox[2]{%
    \newbox\contentbox%
    \newbox\bkgdbox%
    \setbox\contentbox\hbox to \hsize{%
        \vtop{
            \kern\columnsep
            \hbox to \hsize{%
                \kern\columnsep%
                \advance\hsize by -2\columnsep%
                \setlength{\textwidth}{\hsize}%
                \vbox{
                    \parskip=\baselineskip
                    \parindent=0bp
                    #2
                }%
                \kern\columnsep%
            }%
            \kern\columnsep%
        }%
    }%
    \setbox\bkgdbox\vbox{
        \color{#1}
        \hrule width  \wd\contentbox %
               height \ht\contentbox %
               depth  \dp\contentbox
        \color{black}
    }%
    \wd\bkgdbox=0bp%
    \vbox{\hbox to \hsize{\box\bkgdbox\box\contentbox}}%
    \vskip\baselineskip%
}

\date{}

\begin{document}


\title{On the motion of a compressible viscous fluid driven by time periodic inflow/outflow boundary conditions}

\author{Anna Abbatiello\thanks{The research of A.A. is supported by Einstein Foundation, Berlin.} \and Eduard Feireisl\thanks{
The research of E.F. leading to these results has received funding from the
Czech Sciences Foundation (GA\v CR), Grant Agreement
21-02411S. The Institute of Mathematics of the Academy of Sciences of
the Czech Republic is supported by RVO:67985840. The stay of E.F. at TU Berlin is supported by Einstein Foundation, Berlin.
} }


\maketitle

\bigskip

\centerline{Institute of Mathematics, Technische Universit\"{a}t Berlin,}
\centerline{Stra{\ss}e des 17. Juni 136, 10623 Berlin, Germany}
\centerline{anna.abbatiello@tu-berlin.de}

\medskip

\centerline{Institute of Mathematics of the Academy of Sciences of the Czech Republic}
\centerline{\v Zitn\' a 25, CZ-115 67 Praha 1, Czech Republic}
\centerline{feireisl@math.cas.cz}

\begin{abstract}

We consider the barotropic Navier--Stokes system describing the motion of a compressible viscous fluid confined to a bounded domain 
driven by time periodic inflow/outflow boundary conditions. We show that the problem admits a time periodic solution in the class of weak solutions satisfying the energy inequality.

\end{abstract}



\section{Introduction}
\label{I}

Time periodic solutions to \emph{dissipative} dynamical systems result from the balance between the energy dissipation 
and its supply through external driving forces. In many cases, including some simple models in fluid mechanics, the 
problem can be reduced to a system of evolutionary equations with a time periodic right--hand side that represents the influence 
of the outer world, see e.g.  Galdi \cite{Gald}, Maremonti and Padula \cite{MarPad},  the recent survey by Galdi and Kyed \cite{GalKye}, or \cite{Abba}, \cite{AbbMar} and the references therein for the nonlinear viscosity issue. An overwhelming majority of the above cited results concern models of \emph{incompressible} viscous fluids.

Similar problems for compressible fluids have been addressed in \cite{FMPS}, Axmann and Pokorn{\'y} \cite{AXPO}, and in \cite{FeMuNoPo}.    
In all cases the fluid is driven by a time periodic volume force that may be mathematically acceptable but physically 
less relevant. Indeed models of compressible fluids must take into account the mass transport implemented in the system through 
the equation of continuity. The interaction with the outer world is then incorporated in the boundary condition. 

We consider 
a physically realistic scenario when the fluid is driven by general time--periodic inflow/outflow boundary conditions. For the sake of simplicity, we neglect the thermal effects and consider a barotropic fluid, the state of which at a given time instant $t$ and the 
spatial position $x$ is represented by the mass density $\vr = \vr(t,x)$ and the velocity $\vu = \vu(t,x)$. The time evolution of the fluid is determined by the \emph{Navier--Stokes system} of equations:

\begin{equation} \label{I1}
\partial_t \vr + \Div (\vr \vu) = 0,
\end{equation}
\begin{equation} \label{I2}
\partial_t (\vr \vu) + \Div (\vr \vu \otimes \vu) + \Grad p(\vr) = \Div \mathbb{S}(\Grad \vu) + \vr \vc{g}, 
\end{equation}
where the viscous stress $\mathbb{S}$ is given by Newton's rheological law
\begin{equation} \label{I3}
\mathbb{S} (\Grad \vu) = \mu \left( \Grad \vu + \Grad^t \vu - \frac{2}{d} \Div \vu \mathbb{I} \right) + \eta
\Div \vu \mathbb{I},\ \mu > 0, \ \eta \geq 0.
\end{equation}
The pressure $p = p(\vr)$ is an explicitly given function of the density and $\vc{g} = \vc{g}(t,x)$ is a given volume force density. 
Here $\vc{g}$ is allowed to depend on the time but in the real world applications $\vc{g} = \vc{g}(x)$ is just the gravitational force.

The fluid is confined to a bounded domain $\Omega \subset R^d$, on the boundary of which the velocity is determined by general 
inhomogeneous boundary conditions
\begin{equation} \label{I4}
\vu|_{\partial \Omega} = \vuB,\ \vuB = \vuB(t,x). 
\end{equation}
We distinguish the inflow part of the boundary, 
\begin{equation} \label{I5}
\Gamma_{\rm in} = \left\{ (t,x) \ \Big| \ t \in R,\ x \in \partial \Omega,\ \vuB \cdot \vc{n} < 0 \right\}
\end{equation}
and the ouflow part
\begin{equation} \label{I5a}
\Gamma_{\rm out} = \left\{ (t,x) \ \Big| \ t \in R,\ x \in \partial \Omega,\ \vuB \cdot \vc{n} \geq 0 \right\}.
\end{equation}
Both depend on the time $t$ and must be considered as subsets of the time--space cylinder $R \times \partial \Omega$. 
Finally, the density is prescribed on the inflow boundary,
\begin{equation} \label{I6}
\vr = \vrB \ \mbox{on}\ \Gamma_{\rm in}, \ \vrB = \vrB(t,x),\ 0 \leq \vrB < \Ov{\vr}.
\end{equation}
Note that prescribing the density on the whole $\partial \Omega$ the resulting problem is overdetermined and not solvable for any boundary data \eqref{I4} as pointed out in an example given in \cite{CiFeJaPe}. 

We suppose that the data $\vc{g}$, $\vuB$, and $\vrB$ are defined for any $t \in R$, $x \in \Ov{\Omega}$, and that they are \emph{time
periodic} with a period $T > 0$. Equivalently, introducing the flat sphere
\[
\mathcal{T}^1 \equiv [0,T]|_{\{ 0, T \}},
\]
we may suppose $t \in \mathcal{T}^1$.

Our goal is to show that the problem \eqref{I1}--\eqref{I6} admits a time periodic solution with the period $T$. It can be shown by direct manipulation that smooth solutions satisfy the (total) energy balance equation 
\begin{equation} \label{E1}
\begin{split}
&\frac{{\rm d}}{{\rm d}t} \intO{ \left( \frac{1}{2} \vr |\vc{u} - \vuB |^2 + P(\vr) \right) } +
\int_{\partial \Omega} P(\vr) \vuB \cdot \vc{n} \D \sigma_x   \\ &+ \intO{ \mathbb{S} (\Grad (\vc{u} - \vu_B) ) : \Grad (\vc{u} - \vuB) }=- \intO{  \vr (\vc{u} - \vuB) \cdot \Grad \vu_B \cdot (\vc{u} - \vuB) }
 \\
& - \intO{ p(\vr) \Div \vc{u}_B }+ \intO{ \Big( \vr \vc{g} - \vr \partial_t \vuB + \Div \mathbb{S} (\Grad \vu_B)-  \vr \vu_B \cdot \Grad \vu_B \Big) \cdot (\vc{u} -
\vuB) },
\end{split}
\end{equation}
where $P(\vr)$ is the pressure potential satisfying 
\[
P'(\vr) \vr - P(\vr) = p(\vr).
\]
For the anticipated time--periodic motion, the integrals on the right--hand side of \eqref{E1} must be controlled by the dissipation terms on the left--hand side. It turns out that the principal difficulty is to handle possible density concentrations that would render the integral 
\[
\intO{  \vr (\vc{u} - \vuB) \cdot \Grad \vu_B \cdot (\vc{u} - \vuB) }
\]
uncontrollable. Similar problem occurs already in the stationary case, where a suitable remedy is to impose certain 
growth conditions on the pressure $p$. Inspired by \cite{CiFeJaPe}, \cite{CiFeJaPe1}, we consider the \emph{hard sphere pressure}
equation of state coupled with a suitable growth condition,
\begin{equation} \label{R1}
\begin{split}
p \in C^1[0,\Ov{\vr}) \cap C^2(0, \Ov{\vr}),\ p'(\vr) > 0 \ \mbox{for any}\ 0 < \vr < \Ov{\vr},\\
\liminf_{\vr \to \Ov{\vr}-} (\Ov{\vr} - \vr)^\beta p(\vr) > 0 
\end{split}
\end{equation}
for some constants $\beta > 0$,  $\Ov{\vr} > 0$. In particular, the pressure becomes singular as $\vr \nearrow \Ov{\vr}$ and the density is confined to the range $0 \leq \vr < \Ov{\vr}$. The reader may consult Kastler et al. \cite{KVB62CPES}, or 
Kolafa et al. \cite{KLM04AEHS} for the physical background of \eqref{R1}. 

The density being {\it a priori} bounded, the energy balance \eqref{E1} gives rise to {\it a priori} bounds that are sufficient for showing the existence of a time periodic solution. The proof, however, is not completely straightforward, as approximation scheme 
used e.g. in \cite{CiFeJaPe} or \cite{CiFeJaPe1} is not suitable for a time dependent problem. 

As far as the initial--boundary value problem is concerned, there are two approaches available in the literature. Novo \cite{NOVO} and later Girinon \cite{GI} use 
velocity penalization converting the existence proof to the study of a singular limit when the penalization forces the velocity to 
attain the desired value in a small neighborhood of the boundary. This approach is very elegant
and suitable for the initial--boundary value problem as the basic steps of the proof are identical 
with the existing theory for the homogeneous system. The boundary value $\vrB$ of the density is then incorporated in the initial 
data via the method of characteristics. As the initial data of a time--periodic solution are not known {\it a priori}, the application of the penalization method in the present setting is not clear. An alternative approach, similar to \cite{CiFeJaPe}, 
was used by Chang, Jin, Novotn\' y \cite{ChJiNo}: The inhomogeneous boundary conditions are implemented directly in the first step of the approximation process. This requires certain regularity of the boundary $\partial \Omega$ or at least its component $\Gamma_{\rm in}$.

To attack the time--periodic problem, we use the approximation scheme similar to \cite{ChJiNo} at the level of the 
continuity equations, and with a friction type penalization in the momentum balance. The resulting basic approximate problem may be  then solved by a direct method similar to \cite{FeMuNoPo} or employing a fixed-point argument as in \cite{FMPS}. Passing to the limit in the sequence of approximate time--periodic solution requires a non--trivial modification of \cite{FMPS} due to the inhomogeneous boundary conditions. 

The paper is organized as follows. In Section \ref{R}, we state the principal hypotheses and formulate the main result. 
In Section \ref{A}, we introduce the approximation scheme and show the existence of the basic approximate solution in Section \ref{F}. 
The limit in the sequence of approximate solutions is performed in Section \ref{L}. The paper is concluded by a short discussion of possible extensions and further applications of the results in Section \ref{C}.  

\section{Principal hypotheses and the main result}
\label{R}

Before stating our main result, we introduce the basic hypotheses imposed on the data. To avoid technicalities, we suppose that 
$\Omega \subset R^d$ is a bounded domain with a smooth boundary $\partial \Omega$ of class at least $C^3$. In particular, the outer normal vector $\vc{n}(x)$ exists at any $x \in \partial \Omega$, and $\Gamma_{\rm in}$, $\Gamma_{\rm out}$ are well defined. 

\subsection{Boundary velocity decomposition}

If $\widetilde{\vu}_B = \widetilde{\vu}_B(t, x)$ is a boundary velocity field satisfying 
\[
\int_{\Gamma}\widetilde{\vu}_B\cdot \vc{n} \D \sigma_x = 0
\]
for any component $\Gamma \subset \partial \Omega$, then $\widetilde{\vu}_B$ admits an extension inside $\Omega$ in the form
\[
\widetilde{\vu}_B = {\bf curl}_x \vc{w}_B \ \mbox{if}\ d = 3,\ 
\widetilde{\vu}_B = \Grad^\perp w_B, \ \Grad^\perp = (-\partial_{x_2}, \partial_{x_1}) \ \mbox{if}\ d = 2,  
\]
for a certain vector potential $\vc{w}_B$ (or scalar $w_B$), see Galdi \cite[Lemma IX.4.1]{GALN} or Kozono and Yanagisawa 
\cite[Proposition 1]{KozYan}. 

Accordingly, we make the following ansatz for the boundary velocity:
\begin{equation} \label{R2}
\vuB = {\bf curl}[ \vc{w}_B ] + \vc{v}_B \ \mbox{if}\ d = 3,\
\vuB = \Grad^\perp [w] + \vc{v}_B, \ \mbox{if}\ d = 2,
\end{equation}
where $\vc{w}_B = \vc{w}_B(t,x)$, $\vc{v}_B = \vc{v}_B(t,x)$ are smooth, say 
\begin{equation} \label{R2a}
\vc{w}_B \in C^3 (\mathcal{T}^1 \times \Ov{\Omega}; R^d),\ \vc{v}_B 
\in C^2 (\mathcal{T}^1 \times \Ov{\Omega}; R^d),
\end{equation}
and
\begin{equation} \label{R3}
\Ds \vc{v}_B \geq 0,\ \inf_{t,x} \Div \vc{v}_B|_{D} > 0 ,\ D \subset \Omega \ \mbox{open, non--void, where}\ 
\Ds \equiv \frac{1}{2} \left( \Grad + \Grad^t \right).
\end{equation}
where the notation $\Ds \vc{v}_B \geq 0$ means that $\Ds \vc{v}_B$ is a positively semidefinite symmetric matrix.
The component $\vc{v}_B$ provides a stabilizing effect already exploited in \cite{BreFeiNov20}. Alternatively, we may suppose that 
the vector field $\vuB$ is tangential to the boundary, $\vuB \cdot \vc{n} = 0$, and prescribe the total mass 
$M = \intO{ \vr }$.

\subsection{Weak solution}
\label{w}

Before stating the main result, we introduce the concept of weak solution to the time-periodic problem
\eqref{I1}--\eqref{I6}.

\begin{Definition}[Weak solution] \label{Dw1}

We shall say that $[\vr, \vu]$ is a weak time-periodic solution to the problem \eqref{I1}--\eqref{I6} 
if the following holds:

\begin{itemize}

\item {\bf Regularity class.}
\begin{equation} \label{w1}
\begin{split}
0 \leq \vr < \Ov{\vr} \ &\mbox{a.e. in}\ \mathcal{T}^1 \times \Omega,\ \vr \in C_{\rm weak} (\mathcal{T}^1; L^q(\Omega)),\
1 \leq q < \infty,\\
\vr, \ P(\vr) &\in L^1(\Gamma_{\rm out}; |\vuB \cdot \vc{n}| \D \sigma_x \dt );
\end{split}
\end{equation}
\begin{equation} \label{w2}
\vu \in L^2(\mathcal{T}^1; W^{1,2}(\Omega; R^d)),\ (\vu - \vuB) \in L^2(\mathcal{T}^1; W^{1,2}_0(\Omega; R^d));
\end{equation}
\begin{equation} \label{w3}
\vr \vu \in C_{\rm weak}(\mathcal{T}^1; L^2(\Omega; R^d)).
\end{equation}

\item {\bf Equation of continuity.}
\begin{equation} \label{w4}
\intTO{ \Big[ \vr \partial_t \varphi + \vr \vu \cdot \Grad \varphi \Big] } = \int_{\Gamma_{\rm in}} \varphi 
\vrB\vuB \cdot \vc{n}\ \D \sigma_{x} \dt + \int_{\Gamma_{\rm out}} \varphi \vr \vuB \cdot \vc{n}\ \D \sigma_{x} \dt
\end{equation}
for all $\varphi \in C^1(\mathcal{T}^1 \times \Ov{\Omega})$.

\item {\bf Momentum equation.}
\begin{equation} \label{w5}
\intTO{ \Big[ \vr \vu\cdot  \partial_t \bfphi + \vr \vu \otimes \vu : \Grad \bfphi + p(\vr) \Div \bfphi \Big] }
= \intTO{ \Big[ \mathbb{S}(\Grad \vu) : \Grad \bfphi - \vr \vc{g} \cdot \bfphi \Big] }
\end{equation}
for all $\bfphi \in C^1(\mathcal{T}^1; C^1_c ( \Omega; R^d))$.

\item {\bf Energy balance.}
\begin{equation} \label{w6}
\begin{split}
- \intT \partial_t \psi &\intO{ \left[ \frac{1}{2} \vr |\vc{u} - \vuB |^2 + P(\vr) \right] } \dt\\ &+
\int_{\Gamma_{\rm in}} \psi P(\vrB) \vuB \cdot \vc{n} \D \sigma_{x} \dt  +
\int_{\Gamma_{\rm out}} \psi P(\vr) \vuB \cdot \vc{n} \D \sigma_{x} \dt \\&+
\intT \psi \intO{ \mathbb{S} (\Grad (\vc{u} - \vu_B) ) : \Grad (\vc{u} - \vuB) }\dt \\ \leq &-
\intT \psi
\intO{  \vr (\vc{u} - \vuB) \cdot \Grad \vu_B \cdot (\vc{u} - \vuB) } \dt
 - \intT \psi \intO{ p(\vr) \Div \vc{u}_B }\dt \\
&+ \intT \psi \intO{ \Big( \vr \vc{g} - \vr \partial_t \vu_B + \Div \mathbb{S} (\Grad \vu_B)-  \vr \vu_B \cdot \Grad \vu_B \Big) \cdot (\vc{u} -
\vuB) } \dt
\end{split}
\end{equation}
for any $\psi \in C^1(\mathcal{T}^1)$, $\psi \geq 0$.

\end{itemize}

\end{Definition}

Strictly speaking, the available regularity of the density $\vr$ does not guarantee the existence of its trace on $\partial \Omega$. However, the 
velocity $\vu$ being a Sobolev function possesses a well defined trace $\vuB$, while the momentum $\vr \vu$ admits a normal trace 
$\vr \vu \cdot \vc{n}$
in the sense of 
Chen, Torres, and Ziemer \cite{ChToZi}. Accordingly, we may identify $\vr$ with a function in $L^1(\Gamma_{\rm out}; |\vuB \cdot \vc{n}|
\D \sigma_x \dt )$.

\subsection{Main result}

We are ready to state our main result.

\begin{Theorem}[Existence of a time--periodic solution] \label{wT1} 

Let $\Omega \subset R^d$, $d=2,3$ be a bounded domain of class $C^3$.
Let the pressure $p$ satisfy \eqref{R1}, with $\beta \geq 3$. 
Let the boundary velocity $\vuB$ be given through 
\eqref{R2}--\eqref{R3}, and let $\vr_B$ satisfy
\begin{equation}\label{rhob}
\vr_B \in C^1(\mathcal{T}^1 \times \Ov{\Gamma}_{\rm in} ), \ 0 \leq \vr_B < \Ov{\vr}.
\end{equation}
Finally, suppose 
\[
\vc{g} \in L^\infty( \mathcal{T}^1 \times \Omega).
\]

Then the problem \eqref{I1}--\eqref{I6} admits a time-periodic weak solution $[\vr, \vu]$ in the sense of Definition \ref{Dw1}.

\end{Theorem}

For the sake of simplicity, we assume that both the physical domain and the boundary data are regular. This can be certainly relaxed at the expense of additional technicalities in the course of the proof. The rest of the paper is devoted to the proof of Theorem \ref{wT1}.

\section{Approximation scheme}
\label{A}

Before introducing the approximation scheme, we adjust the extension of the boundary velocity $\vuB$.

\subsection{Adjusting the velocity extension}
\label{a}

Let us recall the Korn--Poincar\' e inequality
\begin{equation} \label{a1}
\| \vu - \vuB \|_{W^{1,2}_0(\Omega; R^d)}^2 \leq c_{KP} \intO{ \mathbb{S}(\Grad(\vu - \vuB)) : \Grad (\vu - \vuB) },
\end{equation}
and the Hardy--Sobolev inequality
\begin{equation} \label{a2}
\left( \intO{ \frac{ |\vu - \vuB|^2 }{{\rm dist}^2(x, \partial \Omega)} } \right)^{1/2} \leq c_{HS} \| \vu - \vuB \|_{W^{1,2}_0(\Omega; R^d)},
\end{equation}
where the constants depend only on the geometry of the domain $\Omega$.

Next, we report the following result (see e.g. Galdi \cite[Lemma III.6.1,Lemma III.6.2]{GALN}): \\For each $\omega > 0$, there exists a function
$d_\omega \in C^\infty(\Ov{\Omega})$ enjoying the following properties:
\begin{itemize}
\item
\begin{equation} \label{a3}
|d_\omega | \leq 1,\ d_\omega (x) \equiv 1 \ \mbox{for all}\ x \ \mbox{in an open neighborhood of}\ \partial \Omega;
\end{equation}
\item
\begin{equation} \label{a4}
d_\omega (x) \equiv 0 \ \mbox{whenever}\ {\rm dist}(x, \partial \Omega) > \omega;
\end{equation}
\item
\begin{equation} \label{a5}
|D^\alpha_x d_\omega (x) | \leq c \frac{\omega}{{\rm dist}^{|\alpha|} (x, \partial \Omega)},\ |\alpha| = 1,2,\
x \in \Omega;
\end{equation}

\end{itemize}
where the constant is independent of $\omega$.

In accordance with \eqref{R2} we can choose
\begin{equation}\label{uB}
\vuB = {\bf curl}[ d_\omega \vc{w}_B ] + \vc{v}_B
\end{equation}
for a suitable $\omega > 0$.

It follows from \eqref{a1}, \eqref{a2}, and \eqref{a5} that $\omega > 0$ can be chosen small enough so that
\begin{equation} \label{a6}\begin{split}
- \intO{  \vr (\vc{u} - \vuB) \cdot \Grad \vuB \cdot (\vc{u} - \vuB) }
\leq \frac{1}{4} \intO{ \mathbb{S}(\Grad( \vu - \vuB)) &: \Grad (\vu - \vuB) }\\
\ &\mbox{whenever}\ 0 \leq \vr <\Ov{\vr}.
\end{split}\end{equation}
Indeed, by virtue of \eqref{R3} it holds
$$  \intO{  \vr (\vc{u} - \vuB) \cdot \Grad \vc{v}_B \cdot (\vc{u} - \vuB) }= \intO{ \vr(\vc{u} - \vuB)\otimes(\vc{u} - \vuB): \Ds \vc{v}_B  } \geq 0,$$
then 
\begin{equation*} \begin{split}
&- \intO{  \vr (\vc{u} - \vuB) \cdot \Grad \vuB \cdot (\vc{u} - \vuB) }\\
&= - \intO{  \vr (\vc{u} - \vuB) \cdot \Grad[ {\bf curl}[ d_\omega \vc{w}_B ] ] \cdot (\vc{u} - \vuB) } - \intO{  \vr (\vc{u} - \vuB) \cdot \Grad \vc{v}_B \cdot (\vc{u} - \vuB) } \\
&\leq  \Ov{\vr}\intO{|\vu-\vuB|^2|\Grad[ {\bf curl}(d_\omega \vc{w}_B) ]|} 
\end{split}\end{equation*}
and,  since
\[ |\Grad[ {\bf curl}(d_\omega \vc{w}_B) ]|\leq |\Grad^2 d_\omega| \|\vc{w}_B\|_\infty +  |\Grad d_\omega| \|\Grad\vc{w}_B\|_\infty + |d_\omega| \|\Grad^2\vc{w}_B\|_\infty \leq c \frac{\omega}{{\rm dist}^{2} (x, \partial \Omega)},\]
 employing \eqref{a2} and \eqref{a1}, the assertion follows for a suitable $\omega >0$.
In the remaining part of the paper, we suppose that $\omega > 0$ has been fixed for \eqref{a6} to hold.

\subsection{Approximate equation of continuity}

Similarly to Chang et al. \cite{ChJiNo}, the equation of continuity is approximated as
\begin{equation} \label{A1}
\begin{split}
\partial_t \vr + \Div (\vr \vu) + \ep \vr &= \ep \Del \vr \mbox{ in }(0, T)\times \Omega, \\
\ep \Grad \vr \cdot \vc{n} + (\vrB - \vr) [\vuB \cdot \vc{n}]^- &= 0 \ \mbox{on } (0, T)\times \partial \Omega,
\end{split}
\end{equation}
where $\ep > 0$ is a positive parameter, and $[\vu_B \cdot \vc{n}]^- = \min \left\{ \vu_B \cdot \vc{n} ; 0 \right\}$. Given $\vu$ sufficiently regular and  the data $\vrB, \vuB$, then the existence of a unique solution to \eqref{A1} with initial condition 
\begin{equation*} 
\vr(0)= \vr_0 \mbox{ in } \Omega, \vr_0\geq 0,\end{equation*}
follows from the application of the maximal regularity theory for parabolic initial-boundary value problems with inhomogeneous data.
Such result is obtained for instance in  \cite[Lemma 4.3]{ChJiNo} asking $\vr_0 \in W^{1, 2}(\Omega)$ and employing the maximal regularity theorem by Denk et al. \cite[Theorem 2.1]{DeHiPr}. However for our aims it is enough to first work with the Faedo--Galerkin approximation of problem \eqref{A1} and to this aim we need only to consider the initial density in a finite dimensional subspace of $L^2(\Omega)$. Thus, let us introduce 
$$Y_n = {\rm span} \left\{ {z}_i\ \Big|\ {z}_i \in C^\infty(\Ov{\Omega}),\ i = 1,\dots, n \right\}$$
where ${z}_i$ are orthonormal with respect to the standard scalar product in $L^2$.  We look for the approximate density $\vr \in C^1([0, T]; Y_n)$ such that 
\begin{equation}\label{Ga1}
\begin{split}
\intO{\partial_t\vr \varphi} + \intO{\Div(\vr\vu)\varphi} + \ep\intO{\vr \varphi} +\ep \intO{\Grad\vr\cdot\Grad\varphi}& \\+ \int_{\partial\Omega} (\vrB-\vr) [\vuB \cdot \vc{n}]^- \varphi \D \sigma_{x}=0 \mbox{ for any } \varphi \in Y_n, \ t\in (0, T),&
\end{split}
\end{equation}
Given the initial condition $\vr(0, \cdot)= \vr_0 \in Y_n$, $\vu \in L^\infty((0, T); W^{1, \infty}(\Omega; R^d))$, $\vuB$ and $\vrB$ as in \eqref{uB} and \eqref{rhob} respectively, the existence of the  Faedo-Galerkin approximation follows by the classical theory of ODEs, indeed the problem \eqref{Ga1} is a linear system of ODEs for the unknown $\vr$. Moreover, taking $\varphi=\vr$ it follows
\begin{equation*}
\begin{split}
\frac{1}{2}\frac{{\rm d}}{{\rm dt}}\|\vr\|_{L^2(\Omega)}^2 + \intO{\vr \,&\Div(\vr\vu)} + \varepsilon \|\vr\|_{L^2(\Omega)}^2 + \varepsilon \|\Grad\vr\|_{L^2(\Omega)}^2\\& +\int_{\partial\Omega}\vr (\vrB-\vr) [\vuB \cdot \vc{n}]^- \, \D \sigma_{x} =0.
\end{split}\end{equation*} 
Since $-\int_{\partial\Omega}\vr^2 [\vuB \cdot \vc{n}]^- \D \sigma_{x}\geq 0$, we have 
\begin{equation*}\begin{split}
\frac{1}{2}\frac{{\rm d}}{{\rm dt}}\|\vr\|_{L^2(\Omega)}^2 +\varepsilon \|\vr\|_{L^2(\Omega)}^2 + \varepsilon \|\Grad\vr\|_{L^2(\Omega)}^2 \leq \left|\intO{\vr^2\Div (\vu)}\right|& \\ +  \left|\intO{\vr \Grad\vr\cdot \vu}\right| + \left|\int_{\partial\Omega}\vr \vrB [\vuB \cdot \vc{n}]^- \D \sigma_{x}\right|&
\end{split}\end{equation*} 
and, employing the H\"{o}lder and the Young inequalities, then the embedding of the trace spaces into $W^{1, 2}$ and the smoothness of the data, we conclude through the Gronwall inequality that
\begin{equation}\label{U1}
\sup_{[0, T]} \|\vr\|_{L^2(\Omega)}^2 + \ep \int_0^T \|\vr\|_{W^{1, 2}(\Omega)}^2\dt \leq c
\end{equation}
where $c$ is a positive constant depending only on the data and on the norm of $\vu$. Therefore the Faedo--Galerkin approximation $\vr$, fulfilling \eqref{Ga1}, exists globally in time. 

\subsection{Approximate momentum equation}

The momentum equation is replaced by a Faedo--Galerkin approximation. To this end, consider a finite--dimensional space
\[
X_n = {\rm span} \left\{ \vc{w}_i\ \Big|\ \vc{w}_i \in \DC(\Omega; R^d),\ i = 1,\dots, n \right\}
\]
where $\vc{w}_i$ are orthonormal with respect to the standard scalar product in $L^2$. 

We look for the approximate velocity field in the form
\[
\vu = \vc{v} + \vuB, \ \vc{v} \in C([0,T]; X_n),
\]
where
\begin{equation} \label{A2}
\begin{split}
&\left[ \intO{ (\ep + \vr) \vc{v} \cdot \bfphi } \right]_{t=0}^{t = \tau} \\&=
\int_0^\tau \intO{ \Big[ (\ep + \vr)  \vc{v} \cdot \partial_t \bfphi + \vr \vc{v}\otimes (\vc{v}+\vuB) : \Grad \bfphi
+ p_\delta (\vr) \Div \bfphi - \mathbb{S}(\Grad  \vc{v})   : \Grad \bfphi \Big] }\dt\\
&- \int_0^\tau \intO{ \Big[ \vr \partial_t \vuB + \vr  (\vc{v} + \vu_B)\cdot \Grad \vuB \Big] \cdot \bfphi } \dt +  \int_0^\tau \intO{ \Div \mathbb{S}(\Grad\vu_B)\cdot\bfphi}\dt\\
&- \ep \int_0^\tau \intO{ \Grad \vr \cdot \Grad \vc{v} \cdot \bfphi } \dt 
- \ep \int_0^\tau \intO{ \vr \vc{v} \cdot \bfphi } \dt
\\ &+ \int_0^\tau \intO{ \vr \vc{g} \cdot \bfphi } \dt
- \int_0^\tau \intO{ \Lambda (\vr) \vc{v} \cdot \bfphi  } \dt
\end{split}
\end{equation}
for any $\bfphi \in C^1([0,T]; X_n)$.
Here, we have set
\begin{equation} \label{A3}
p_\delta (\vr) = \left\{ \begin{array}{l} p(\vr) + \delta \vr^2 \ \mbox{for}\ 0 \leq \vr \leq \Ov{\vr} - \delta \\ \\
(\vr - \Ov{\vr} + \delta)^\Gamma + p (\Ov{\vr} - \delta) + p' (\Ov{\vr} - \delta) (\vr - \Ov{\vr} + \delta)
 + \delta \vr^2,\\ \  \vr > \Ov{\vr} - \delta
\end{array} \right.
\end{equation}
for $\Gamma > 2$ large enough, and
\begin{equation} \label{A4}
\Lambda(\vr) = \Ov{\Lambda} \left[ \vr - \frac{3}{2} \Ov{\vr} \right]^+, 
\left[ \vr - \frac{3}{2} \Ov{\vr} \right]^+ = \max \left\{ \vr - \frac{3}{2} \Ov{\vr} ; 0 \right\}, 
\end{equation}
with $\Ov{\Lambda}$ -- a positive constant to be fixed below. It is easy to observe (see e.g. \cite[Chapter 7]{EF70}) that 
for given $\vr$, the problem \eqref{A2} represents a (nonlinear) system of ODE's for the unknown $\vc{v}$.

\subsection{Approximate energy inequality}

The energy balance for the approximate scheme can be obtained
by using $\vc{v} = \vu - \vuB$ as a test function in \eqref{A2}:

\begin{equation} \label{F1}
\begin{split}
\frac{{\rm d}}{{\rm d}t}  &\intO{ \left[ \frac{1}{2} (\vr + \ep)  |\vc{u} - \vuB |^2 + P_\delta(\vr) \right] }  + \ep \intO{ \left[ \vr P'_\delta (\vr)  + 
P''_\delta (\vr) |\Grad \vr|^2 \right] } \\&+
\int_{\Gamma_{\rm in}(t)} P_\delta (\vrB) \vuB \cdot \vc{n} \D \sigma_{x}   +
\int_{\Gamma_{\rm out}(t) } P_\delta (\vr) \vuB \cdot \vc{n} \D \sigma_{x}  \\&+
\intO{ \mathbb{S} (\Grad (\vc{u} - \vu_B) ) : \Grad (\vc{u} - \vuB) } + \intO{ \Lambda(\vr) |\vu - \vuB|^2 } + 
\frac{\ep}{2} \intO{ \vr |\vu - \vuB|^2 } \\ &- 
 \int_{\Gamma_{\rm in}(t)} \Big[ P_\delta (\vrB) - P'_\delta(\vr) (\vr_B - \vr) - P_\delta (\vr) \Big] \vuB \cdot \vc{n} \ \D \sigma_{x}
\\  = &-
\intO{  \vr (\vc{u} - \vuB) \cdot \Grad \vu_B \cdot (\vc{u} - \vuB) } 
 - \intO{ p_\delta (\vr) \Div \vc{u}_B } \\
&+ \intO{ \Big( \vr \vc{g} - \vr \partial_t \vu_B + \Div \mathbb{S} (\Grad \vu_B)-  \vr \vu_B \cdot \Grad \vu_B \Big) \cdot (\vc{u} -
\vuB) } 
\end{split}
\end{equation}

Now, there are two crucial observations: 

\begin{enumerate} 

\item 
\begin{equation} \label{F1-1}
- \intO{ p_\delta (\vr) \Div \vc{u}_B }\dt = - \intO{ p_\delta (\vr) \Div \vc{v}_B }\dt \leq 0
\end{equation}
in view of \eqref{R3}; 

\item 
\begin{equation} \label{F2}
\begin{split}
-
&\intO{  \vr (\vc{u} - \vuB) \cdot \Grad \vu_B \cdot (\vc{u} - \vuB) } \leq \\  
&\frac{1}{2} \left[ \intO{ \mathbb{S} (\Grad (\vc{u} - \vu_B) ) : \Grad (\vc{u} - \vuB) } + \intO{ \Lambda(\vr) |\vu - \vuB|^2 }    \right]
\end{split}
\end{equation}
in view of \eqref{a6}, \eqref{A4}, if $\Ov{\Lambda}$ in \eqref{A4} is chosen large enough. 

\end{enumerate}

Therefore, adopting the standard approach we computed  $\vr$ in terms of $\vu=\vc{v}+\vuB,$ $\vc{v}\in C(0, T; X_n)$, and then we can calculate $\vc{v}$, solving \eqref{A2} with given initial condition $\vc{v}(0, \cdot)= \vc{v}_0\in X_n$  applying a fixed point argument following the steps in \cite[Section 4.3.2]{ChJiNo}. Moreover, with the help of the Gronwall inequality one can claim that $\vc{v}$ enjoys 
\begin{equation}\label{U2}
\sup_{t\in [0, T]} \intO{ \left[ \frac{1}{2} (\vr + \ep)  |\vc{u} - \vuB |^2 + P_\delta(\vr) \right] } \leq c
\end{equation}
where $c$ is a positive constant depending only on the data. Therefore the Faedo-Galerkin approximations $\vc{v}$ is globally well defined. 

\section{First approximation level}
\label{F}

Keeping $n > 0$, $\ep > 0$, and $\delta > 0$ fixed, we aim to show the existence of a time--periodic solution to the approximate system.
 
\subsection{Time-periodic approximations}

The existence of time-periodic approximations will follow from an application of the Brouwer fixed-point theorem to the following mapping 
\[ \begin{split}
\mathcal{F} : Y_n \times X_n &\to Y_n \times X_n,\\
\mathcal{F}[\vr(0), \vc{v}(0)]& = [\vr(T), \vc{v}(T)],
\end{split}  \]
where $(\vr(t), \vc{v}(t))$ are the solutions to \eqref{Ga1}, \eqref{A2} with initial data $\vr(0), \vc{v}(0)$ respectively. 
First, note that $\mathcal{F}$ is well-defined thanks to the existence and uniqueness of the approximating solutions $(\vr(t), \vc{v}(t))$. Indeed the existence theory is discussed in the previous section and however the proof of the uniqueness follows the same lines than the proof of the continuity of $\mathcal{F}$ performed in this section.

Now, let us observe that from \eqref{F1}, employing \eqref{F1-1}, \eqref{F2}, and the facts that $P_\delta$ is a convex function and $P_\delta (\vr) \vuB \cdot \vc{n}$ is non-negative on $\Gamma_{\rm out}$,  it follows 
\begin{equation} \label{}
\begin{split}
&\frac{{\rm d}}{{\rm d}t} \intO{ \left[ \frac{1}{2} (\vr + \ep)  |\vc{v} |^2 + P_\delta(\vr) \right] }  + \ep \intO{ \vr P'_\delta (\vr) }   
+
\frac{1}{2}\intO{ \mathbb{S} (\Grad \vc{v}) : \Grad \vc{v}  } \\&+ 
\frac{\ep}{2} \intO{ \vr |\vc{v}|^2 } 
\leq  \int_{\Gamma_{\rm in}(t) } |P_\delta (\vrB) \vuB \cdot \vc{n}| \D \sigma_{x} \\&+  \intO{ \Big[ \vr (|\vc{g}| + | \partial_t \vu_B| +  | \vu_B \cdot \Grad \vu_B|) + |\Div \mathbb{S} (\Grad \vu_B)| \Big]  |\vc{v}| };
\end{split}
\end{equation}
then using the H\"{o}lder and Young inequalities, \eqref{a1}, the integrability property $\vr\in C(0, T; L^2(\Omega))$ and the fact that $ P'_\delta (\vr) \vr = P_\delta(\vr) + p_\delta(\vr)\geq P_\delta(\vr)$, we obtain
\begin{equation} \label{AF1}
\begin{split}
&\frac{{\rm d}}{{\rm d}t} \intO{ \left[ \frac{1}{2} (\vr + \ep)  |\vc{v} |^2 + P_\delta(\vr) \right] } + \ep \intO{ \left[ \frac{1}{2} (\vr + \ep)  |\vc{v} |^2 + P_\delta(\vr) \right] }
\leq c
\end{split}
\end{equation}
where $c$ is a positive constant depending on the data. Thus, 
after the integration in time we get the existence of ${E}={E}(n, \delta, \ep)$  such that
$$  \intO{ \left[ \frac{1}{2} (\vr + \ep)  |\vc{v} |^2 + P_\delta(\vr) \right]\!\! \!}\ (0)\leq E \ \Longrightarrow  \intO{ \left[ \frac{1}{2} (\vr + \ep)  |\vc{v} |^2 + P_\delta(\vr) \right]\!\! \!}\ (T)\leq E;$$
next employing $\vr\geq 0$, the equivalence of the norms on $X_n$, the definition \eqref{A3} and the equivalence of the norms on $Y_n$, we derive that 
$$ \|\vc{v}(0)\|_{X_n}\leq {E} \ \Longrightarrow  \ \|\vc{v}(T)\|_{X_n}\leq {E}$$
and 
$$\|\vr(0)\|_{Y_n}\leq E \ \Longrightarrow \ \|\vr(T)\|_{Y_n}\leq E.$$
Consequently $\mathcal{F}$ maps the following set  into itself
$$\mathcal{K}[{E}]:= \left\{ (r,  \vc{w}) \in Y_n \times X_n : \|r\|_{Y_n }\leq E,\ \|\vc{w}\|_{X_n}\leq {E} \right\}. $$

Finally, we need to prove the continuity of $\mathcal{F}$. Let $\vr_1, \vr_2$ be solutions to \eqref{Ga1} with  initial data ${\vr_1}(0), {\vr_2}(0)$ and,  let $\vc{v}_1, \vc{v}_2$ be solutions to \eqref{A2} with initial data $\vc{v}_1(0), \vc{v}_2(0)$. Let us take the difference among the Galerkin formulations  \eqref{Ga1} for $\vr_1, \vr_2$:
\begin{equation*}
\begin{split}
& \intO{ \partial_t(\vr_1-\vr_2) \varphi } +\intO{\varphi\Div(\vr_1\vu_1-\vr_2\vu_2)} + \varepsilon \intO{(\vr_1-\vr_2)\varphi} \\&+\varepsilon \intO{\Grad(\vr_1-\vr_2)\cdot\Grad\varphi} - 
\int_{\partial\Omega} \varphi (\vr_1 -\vr_2) [\vu_B \cdot \vc{n}]^- \ \D  \sigma_{x} =0,
\end{split}
\end{equation*}
then choose $\varphi =\vr_1-\vr_2$
\begin{equation*}
\begin{split}
&\Dt \intO{\frac{1}{2}|\vr_1-\vr_2|^2} +\intO{(\vr_1-\vr_2)\Div(\vr_1\vu_1-\vr_2\vu_2) } + \varepsilon \intO{|\vr_1-\vr_2|^2} \\&+\varepsilon \intO{|\Grad(\vr_1-\vr_2)|^2} - 
\int_{\partial\Omega} |\vr_1 -\vr_2|^2 [\vu_B \cdot \vc{n}]^- \ \D  \sigma_{x} =0,
\end{split}
\end{equation*}
since $ - \int_{\partial\Omega} |\vr_1 -\vr_2|^2 [\vu_B \cdot \vc{n}]^- \ \D  \sigma_{x}\geq 0$ and using the uniform estimates  \eqref{U1}, \eqref{U2}, it follows
\begin{equation}
\Dt \intO{\frac{1}2 |\vr_1-\vr_2|^2} +  \intO{|\vr_1-\vr_2|^2} \aleq 1.
\end{equation}
Thus, there exists $\eta=\eta(n, \varepsilon, \delta)>0$ such that
\begin{equation}\label{c1}
 \|\vr_1(0)-\vr_2(0)\|_{L^2(\Omega)}\leq \eta \ \Longrightarrow \  \|\vr_1(T)-\vr_2(T)\|_{L^2(\Omega)}\leq \eta. \end{equation}
Now, analogously let us consider the Galerkin formulations \eqref{A2} for $\vc{v}_1, \vc{v}_2$ and take the difference
\begin{equation*}\begin{split}
\varepsilon\intO{\partial_t(\vc{v}_1-\vc{v}_2) \cdot \bfphi} + \intO{(\mathbb{S}(\Grad  \vc{v}_1) -\mathbb{S}(\Grad  \vc{v}_2)):\Grad \bfphi} = - \intO{(\vr_1\partial_t\vc{v}_1-\vr_2\partial_t\vc{v}_2) \cdot \bfphi}&\\
- \intO{(\vr_1\vu_1\cdot\Grad\vu_1 - \vr_2\vu_2\cdot\Grad\vu_2)\cdot\bfphi} + \varepsilon \intO{(\Grad\vr_1-\Grad\vr_2)\cdot\Grad\bfphi \cdot \vc{v}_1}&\\
 + \varepsilon \intO{\Grad\vr_2\cdot\Grad\bfphi\cdot(\vc{v}_1-\vc{v}_2)} + \intO{ (p_\delta (\vr_1)- p_\delta (\vr_2))\Div\bfphi} &\\
- \intO{(\vr_1 - \vr_2) \partial_t\vuB \cdot \bfphi} + \intO{(\vr_1-\vr_2)\vc{g} \cdot \bfphi} -\intO{ (\Lambda (\vr_1) \vc{v}_1 - \Lambda (\vr_2) \vc{v}_2) \cdot \bfphi  },&
\end{split}\end{equation*}
choose $\bfphi = \vc{v}_1 - \vc{v}_2$ and then arguing as in \eqref{AF1} in view of the uniform estimates obtained after the integration in time of \eqref{F1}, we deduce 
\begin{equation*}
\frac{\varepsilon}{2} \Dt \intO{|\vc{v}_1-\vc{v}_2|^2} + \intO{|\vc{v}_1-\vc{v}_2|^2} \aleq 1
\end{equation*}
as a consequence it follows the existence of $\eta=\eta(n, \delta, \varepsilon)$ such that
\begin{equation}\label{c2}
\|\vc{v}_1(0)-\vc{v}_2(0)\|_{X_n} \leq \eta \ \Longrightarrow \ \|\vc{v}_1(T)-\vc{v}_2(T)\|_{X_n} \leq \eta.
\end{equation}
Therefore \eqref{c1} and \eqref{c2} prove that $\mathcal{F}$ is continuous. 

Consequently, we conclude that there exists a time--periodic 
solution $[\vr, \vu]$ at the first approximation level. 
Moreover, it holds the energy balance 

\begin{equation} \label{F3}
\begin{split}
- \int_{\mathcal{T}^1} \partial_t \psi  &\intO{ \left[ \frac{1}{2} (\vr + \ep)  |\vc{u} - \vuB |^2 + P_\delta(\vr) \right] } \dt + 
\ep \int_{\mathcal{T}^1} \psi \intO{ \Big[ \vr P'_\delta (\vr)  + 
P''_\delta (\vr) |\Grad \vr|^2 \Big] } \dt
  \\&+
\int_{\Gamma_{\rm in} } \psi P_\delta (\vrB) \vuB \cdot \vc{n} \D \sigma_{x} \dt   +
\int_{\Gamma_{\rm out}  } \psi P_\delta (\vr) \vuB \cdot \vc{n} \D \sigma_{x} \dt  \\&+
\int_{\mathcal{T}^1} \psi \intO{ \mathbb{S} (\Grad (\vc{u} - \vu_B) ) : \Grad (\vc{u} - \vuB) } \dt \\&+ \int_{\mathcal{T}^1} 
\psi \intO{ \Lambda(\vr) |\vu - \vuB|^2 }\dt + 
\frac{\ep}{2} \int_{\mathcal{T}^1} \psi \intO{ \vr |\vu - \vuB|^2 }\dt \\ &- 
 \int_{\Gamma_{\rm in} } \psi \Big[ P_\delta (\vrB) - P'_\delta(\vr) (\vr_B - \vr) - P_\delta (\vr) \Big] \vuB \cdot \vc{n} \ \D \sigma_{x} \dt 
\\  = &-
\int_{\mathcal{T}^1} \psi \intO{  \vr (\vc{u} - \vuB) \cdot \Grad \vu_B \cdot (\vc{u} - \vuB) } \dt 
 - \int_{\mathcal{T}^1} \psi  \intO{ p_\delta (\vr) \Div \vc{u}_B } \dt \\
&+ \int_{\mathcal{T}^1} \psi \intO{ \Big( \vr \vc{g} - \vr \partial_t \vu_B + \Div \mathbb{S} (\Grad \vu_B)-  \vr \vu_B \cdot \Grad \vu_B \Big) \cdot (\vc{u} -
\vuB) } \dt
\end{split}
\end{equation}
for any $\psi \in C^1(\mathcal{T}^1)$.\\

Note that instead of the fixed point argument of \cite{FMPS}, we could have used the direct argument of \cite{FeMuNoPo} to prove the existence of a time--periodic approximation.

\section{Asymptotic limit} 
\label{L} 

There are three levels of limits to be performed in this order: $n \to \infty$, $\ep \to 0$, and $\delta \to 0$. 
The limit $n \to \infty$ is nowadays well understood and can be carried out in a way similar to Chang et al. \cite{ChJiNo} 
or \cite{FMPS}. The limits $\ep \to 0$ and $\delta \to 0$ are quite similar, except handling the pressure perturbation in the latter case. We therefore focus on the most difficult last limit $\delta \to 0$. 
Accordingly, we suppose there is a family $\{ \vrd, \vud \}_{\delta > 0}$
of time periodic solutions satisfying: 

\begin{itemize}

\item 
\begin{equation} \label{L1}
\begin{split}
0 &\leq \vrd \ \mbox{ a.e. in}\ \mathcal{T}^1 \times \Omega,\ \vrd \in C_{\rm weak} (\mathcal{T}^1; L^\Gamma(\Omega)),\
 \\ \vrd,\ 
P_\delta (\vrd) &\in L^1(\Gamma_{\rm out}; |\vuB \cdot \vc{n}| \D \sigma \dt ),\\ 
\vud &\in L^2(\mathcal{T}^1; W^{1,2}(\Omega; R^d)),\ (\vud- \vuB) \in L^2(\mathcal{T}^1; W^{1,2}_0(\Omega; R^d)), \\
\vrd \vud &\in C_{\rm weak}(\mathcal{T}^1; L^{\frac{2 \Gamma}{\Gamma + 1}}(\Omega; R^d)); 
\end{split}
\end{equation}

\item
\begin{equation} \label{L2}
\intTO{ \Big[ \vrd \partial_t \varphi + \vrd \vud \cdot \Grad \varphi \Big] } = \int_{\Gamma_{\rm in}} \varphi
\vrB \vuB \cdot \vc{n}\ \D \sigma_{x} \dt +  \int_{\Gamma_{\rm out}} \varphi \vrd \vuB \cdot \vc{n}\ \D \sigma_{x} \dt
\end{equation}
for all $\varphi \in C^1(\mathcal{T}^1 \times \Ov{\Omega})$;

\item
\begin{equation} \label{L3}
\begin{split}
&\intTO{ \Big[ \vrd \vud\cdot \partial_t \bfphi + \vrd \vud \otimes \vud : \Grad \bfphi + p_\delta (\vrd) \Div \bfphi \Big] }\\
&= \intTO{ \Big[ \mathbb{S}(\Grad \vud) : \Grad \bfphi - \Lambda (\vrd) \vud \cdot \bfphi + \vrd \vc{g} \cdot \bfphi \Big] }
\end{split}
\end{equation}
for all $\bfphi \in C^1_c (\mathcal{T}^1 \times \Omega; R^d)$;

\item 
\begin{equation} \label{L4}
\begin{split}
- &\intT  \partial_t \psi \intO{ \left[ \frac{1}{2} \vrd |\vud - \vuB |^2 + P_\delta(\vr) \right] } \dt\\ &+
\int_{\Gamma_{\rm in}} \psi P_\delta (\vrB) \vuB \cdot \vc{n} \D \sigma_{x} \dt  +
\int_{\Gamma_{\rm out}} \psi P_\delta (\vrd) \vuB \cdot \vc{n} \D \sigma_{x} \dt \\&+
\intT \psi \intO{ \mathbb{S} (\Grad (\vud - \vuB) ) : \Grad (\vud - \vuB) }\dt  
+ \intT \psi \intO{ \Lambda(\vrd) |\vud - \vuB|^2 } \dt \\ &\leq -
\intT \psi
\intO{  \vrd (\vud - \vuB) \cdot \Grad \vuB \cdot (\vud - \vuB) } \dt
 - \intT \psi \intO{ p_\delta (\vrd) \Div \vuB }\dt \\
&+ \intT \psi \intO{ \Big( \vrd \vc{g} - \vrd \partial_t \vuB + \Div \mathbb{S} (\Grad \vuB)-  \vrd \vuB \cdot \Grad \vuB \Big) \cdot (\vud -
\vuB) } \dt
\end{split}
\end{equation}
for any $\psi \in C^1(\mathcal{T}^1)$, $\psi \geq 0$.

\end{itemize}

\subsection{Uniform bounds}

Our goal is to derive uniform bounds for the approximate solution $\{ \vrd, \vud \}_{\delta > 0}$ independent of $\delta \to 0$.

\subsubsection{Energy bounds} 

The choice $\psi = 1$ in \eqref{L4} yields
\[
\begin{split}
\intT  \int_D p_\delta (\vrd) \dx \dt &+ 
 \int_{\Gamma_{\rm out}}  P_\delta (\vrd) \vuB \cdot \vc{n} \D \sigma_{t,x} \\&+
\frac{1}{2} \intT \intO{ \mathbb{S} (\Grad (\vud - \vuB) ) : \Grad (\vud - \vuB) }\dt  
+ \frac{1}{2} \intT  \intO{ \Lambda(\vrd) |\vud - \vuB|^2 } \dt \\ \leq &-
 \int_{\Gamma_{\rm in}}  P_\delta (\vrB) \vuB \cdot \vc{n} \D \sigma_{t,x} 
  \\
&+ \intT  \intO{ \Big( \vrd \vc{g} - \vrd \partial_t \vuB + \Div \mathbb{S} (\Grad \vuB)-  \vrd \vuB \cdot \Grad \vuB \Big) \cdot (\vud -
\vuB) } \dt,
\end{split}
\]
where we have used \eqref{F2} and hypothesis \eqref{R3}. Moreover, by virtue of hypothesis \eqref{I6}, 
\[
-\int_{\Gamma_{\rm in}}  P_\delta (\vrB) \vuB \cdot \vc{n} \D \sigma_{t,x} \aleq 1.
\]
Thus using H\" older inequality we may infer that 
\begin{equation} \label{L5}
\begin{split}
\intT  \int_D &p_\delta (\vrd) \dx \dt + 
\int_{\Gamma_{\rm out}}  P_\delta (\vrd) \vuB \cdot \vc{n} \D \sigma_{t,x} \\&+
\frac{1}{4} \intT \intO{ \mathbb{S} (\Grad (\vud - \vu_b) ) : \Grad (\vud - \vuB) }\dt  
+ \frac{1}{4} \intT  \intO{ \Lambda(\vrd) |\vud - \vuB|^2 } \dt \\ &\leq  C \left( 1 + \intTO{ \vrd } \right).
 \end{split}
\end{equation}
At this stage, we do not control the total mass (the integral on the right--hand side) and more elaborated pressure estimates are needed.

\subsubsection{Pressure estimates, I}

At this stage, we need a suitable inverse of the divergence operator. We make use of the construction due to Bogovskii 
\cite{BOG} and introduce the operator $\mathcal{B}$ enjoying the following properties, see Gei{\ss}ert, Heck, and Hieber 
\cite{GEHEHI}: 
\[
\begin{split}
\mathcal{B} : L^q_0 (\Omega) &\equiv \left\{ f \in L^q(\Omega) \ \Big| \ \intO{ f } = 0 \right\}
\to W^{1,q}_0(\Omega, R^d),\ 1 < q < \infty,\\
\Div \mathcal{B}[f] &= f .
\end{split}
\]

Consider a function 
\[
\phi = \phi(x),\ \phi \in C^\infty(\Ov{\Omega}),\ \phi|_{\Omega \setminus D} = 1,\ \intO{\phi} = 0.
\]
Use 
\[
\mathcal{B}[\phi], \ \mbox{where}\ \mathcal{B} \ \mbox{is the Bogovskii operator,}
\]
as a test function in the momentum balance \eqref{L3}:
\begin{equation} \label{L6}
\begin{split}
\intT \int_{\Omega \setminus D } p_\delta (\vrd)  \dx \dt  
&+ \intT \int_D p_\delta (\vrd) \phi \dx \dt
 \\
&= \intTO{ \Big[ \mathbb{S}(\Grad \vud) : \Grad \mathcal{B}[\phi] - \Lambda (\vrd) \vud \cdot \mathcal{B}[\phi] + \vrd \vc{g} \cdot \mathcal{B}[\phi] \Big] } \\
&- \intTO{  \vrd \vud \otimes \vud : \Grad \mathcal{B}[\phi] }.
\end{split}
\end{equation}
Multiplying \eqref{L6} on a small positive constant and adding the resulting expression to \eqref{L5} we deduce the estimate
\begin{equation} \label{L7}
\begin{split}
\intTO{ p_\delta& (\vrd) } + 
\int_{\Gamma_{\rm out}}  P_\delta (\vrd) \vuB \cdot \vc{n} \D \sigma_{x} \dt \\&+
\intT \intO{ \mathbb{S} (\Grad (\vud - \vuB) ) : \Grad (\vud - \vuB) }\dt  
+ \intT  \intO{ \Lambda(\vrd) |\vud - \vuB|^2 } \dt \\ &\aleq  \left( 1 + \intTO{ \vrd } \right).
 \end{split}
\end{equation}

Now, since the pressure becomes singular as $\vr \nearrow \Ov{\vr}$ and enjoys the property in \eqref{R1} with $\beta \geq 3$, passing to a suitable subsequence $\delta \to 0$ as the case may be, we may suppose 
\[
p'(\Ov{\vr} - \delta) \nearrow \infty \ \mbox{as}\ \delta \to 0. 
\]
Consequently, 
a short inspection of \eqref{A3} yields 
\begin{equation} \label{L8}
\frac{ \intTO{ 1_{\vrd \geq \Ov{\vr}  } \vrd } }{\intTO{ p_\delta (\vrd) } } \to 0 \ \mbox{as}\ \delta \to 0.
\end{equation}
Thus \eqref{L7} gives rise to  uniform bounds 
\begin{equation} \label{L9}
\begin{split}
\intTO{ p_\delta (\vrd) } + &
\int_{\Gamma_{\rm out}}  P_\delta (\vrd) \vuB \cdot \vc{n} \D \sigma_{t,x} \\&+
\intT \intO{ \mathbb{S} (\Grad (\vud - \vuB) ) : \Grad (\vud - \vuB) }\dt  \\&
+ \intT  \intO{ \Lambda(\vrd) |\vud - \vuB|^2 } \dt \aleq 1 \ \mbox{as}\ \delta \to  0, 
 \end{split}
\end{equation}
and 
\begin{equation} \label{L10} 
\intTO{ 1_{\vrd \geq \Ov{\vr} } \vrd^\gamma } \to 0 \ \mbox{as}\ \delta \to 0 
\ \mbox{for any}\ 1 \leq \gamma < \Gamma.
\end{equation}

\subsection{Limit $\delta \to 0$}

Our ultimate goal is to perform the limit in the sequence of approximate solutions $\{ \vr_\delta, \vu_\delta \}_{\delta > 0}$. 
To begin, observe that the uniform bounds \eqref{L9}, \eqref{L10}, together with the energy inequality \eqref{L4}, imply also the standard energy estimates 
\begin{equation} \label{L11}
\sup_{t \in \mathcal{T}^1} \left[ \left\| \vrd |\vud|^2 (t, \cdot) \right\|_{L^1(\Omega)} + \left\| P_\delta (\vrd) 
(t, \cdot) \right\|_{L^1(\Omega)} \right] 
\aleq 1.
\end{equation} 
Passing to a suitable subsequence if necessary we may therefore assume that 
\begin{equation} \label{L12}
\begin{split}
\vrd &\to \vr \ \mbox{in}\ C_{\rm weak}(\mathcal{T}^1; L^\Gamma (\Omega)),\ \mbox{where}\ 0 \leq \vr \leq \Ov{\vr}, \\
\vrd &\to \vr \ \mbox{weakly in}\ L^\Gamma (\Gamma_{\rm out}; |\vuB \cdot \vc{n}| \D \sigma_x \dt ),\ 0 \leq \vr \leq \Ov{\vr},\\ 
\int_{\Gamma_{\rm out}} P(\vr) \vuB \cdot \vc{n} \D \sigma_{x} \dt &\leq \liminf_{\delta \to 0} 
\int_{\Gamma_{\rm out}} P_\delta (\vrd) \vuB \cdot \vc{n} \D \sigma_{x} \dt,\\ 
\vud &\to \vu \ \mbox{weakly in}\ L^2(\mathcal{T}^1; W^{1,2}(\Omega; R^d)),\ \vu - \vuB \in L^2(\mathcal{T}^1; W^{1,2}_0(\Omega; R^d)),\\
\vrd \vud &\to \vr \vu \ \mbox{in}\ C_{\rm weak}(\mathcal{T}^1; L^{\frac{2 \Gamma}{\Gamma + 1}}(\Omega; R^d)). 
\end{split}
\end{equation}

In addition, as a consequence of \eqref{L10}, we get 
\begin{equation} \label{L13}
\Lambda (\vrd) \vud \to 0 \ \mbox{in}\ L^q(\mathcal{T}^1 \times \Omega; R^d) \ \mbox{for some}\ q > 1.
\end{equation}

\subsubsection{Renormalized equation of continuity}

As we have shown above, the limit functions $[\vr, \vu]$ satisfy the equation of continuity \eqref{w4}. Unfortunately, this is not enough to perform the last step of the convergence proof -- pointwise convergence of approximate densities. To this end, we need a renormalized version of \eqref{w4}, specifically, 
\begin{equation} \label{L15}
\intTO{ \left[ \beta (\vr) \partial_t \varphi + \beta (\vr) \vu \cdot \Grad \varphi + 
\Big(\beta (\vr) - \beta'(\vr) \vr \Big) \Div \vu \varphi \right]  } = \int_{\Gamma_{\rm in}} 
\varphi \beta(\vr_b) \vuB \cdot \vc{n} \ \D \sigma_{x} \dt
\end{equation} 
for any $\varphi \in C^1_c (\mathcal{T}^1 \times (\Omega \cup \Gamma_{\rm in}))$, $\beta \in BC[0, \infty)$, 
$\beta' \in C_c[0, \infty)$. 

As the limit density is uniformly bounded, relation \eqref{L15} can be shown by the original regularizing argument of 
DiPerna and Lions \cite{DL}. The only problem here is to accommodate the inhomogeneous boundary conditions. Chang et al. 
\cite[Lemma 3.1]{ChJiNo} show \eqref{L15} in the case of time indepedent boundary data $\vrB$, $\vuB$. The proof in the time--dependent case is similar and may be performed via several steps: 
\begin{enumerate}
\item Consider a normal vector field $[\vuB \cdot \vc{n}]^- \vc{n}$ defined on $\partial \Omega$.
\item As $\partial \Omega$ is smooth, it admits an open neighborhood $\mathcal{U}$, $\partial \Omega \subset \mathcal{U}$ 
such that for any $x \in \mathcal{U}$ there is a unique $x_b (x) \in \partial \Omega$ -- the boundary point nearest to $x$.
\item Consider a vector field $\vc{u}_\infty$ defined as 
\[
\vc{u}_{\infty}(t,x) = [\vuB (t, x_b (x)) \cdot \vc{n} (x_b(x))]^- \vc{n} (x_b(x)) \ \mbox{for any}\ x \in \mathcal{U}.
\]
\item Use the method of characteristics to find a solution $\vr_\infty$ of the transport equation
\[
\partial_t \vr_\infty + \Div (\vr_\infty \vc{u}_\infty) = 0 \ \mbox{in}\ \mathcal{T}^1 \times \mathcal{U} \cap (R^d \setminus \Ov{\Omega}),\ 
\vr_\infty|_{\partial \Omega} = \vrB.
\] 

\item As 
\[
\vr_\infty = \vrB,\ \vu_\infty \cdot \vc{n} = \vuB \cdot \vc{n} \ \mbox{on}\ \Gamma_{\rm in}, 
\]
there is an open set $Q \subset \mathcal{T}^1 \times R^d$ such that 
\[
(\mathcal{T}^1 \times (\Omega \cup \Gamma_{\rm in})) \subset Q,
\]
and the functions 
\[
\vr_Q = \left\{ \begin{array}{l} \vr \ \mbox{in} \ (\mathcal{T}^1 \times (\Omega \cup \Gamma_{\rm in})) \\ \\
\vr_\infty \ \mbox{in}\ Q \setminus (\mathcal{T}^1 \times (\Omega \cup \Gamma_{\rm in}))  \end{array} \right. , \ 
\vu_Q = \left\{ \begin{array}{l} \vu \ \mbox{in} \ (\mathcal{T}^1 \times (\Omega \cup \Gamma_{\rm in})) \\ \\
\vu_\infty \ \mbox{in}\ Q \setminus (\mathcal{T}^1 \times (\Omega \cup \Gamma_{\rm in}))  \end{array} \right.
\]
represent a weak solution of the equation of continuity in $Q$.

\item Apply regularizing kernels (both in time and space) to $\vr_Q$, $\vu_Q$ and use the approach to DiPerna and Lions 
\cite{DL} to deduce that $\vr_Q$, $\vu_Q$ is a renormalized solution in $Q$. 

\item As $\vr_Q$, $\vu_Q$ are smooth outside $\mathcal{T}^1 \times \Ov{\Omega}$, we deduce \eqref{L15}.

\end{enumerate} 

\begin{Remark} \label{renR}

Validity of \eqref{L15} can be extended to test functions in the class 
\[
\varphi \in W^{1, \infty}(\mathcal{T}^1 \times \Ov{\Omega}), \ \varphi|_{\Gamma_{\rm out}} = 0
\]
by a density argument.

\end{Remark}

\subsubsection{Pressure estimates II, compactness of the density}

Our ultimate goal is to show equi--integrability of the pressure sequence $\{ p_\delta (\vrd) \}_{\delta > 0}$, and strong (pointwise) 
convergence of $\{ \vrd \}_{\delta > 0}$. As the $L^1-$bound on the pressure has already been established in \eqref{L9}, equi--integrability of $\{ p_\delta (\vrd) \}_{\delta > 0}$ can be shown similarly to \cite[Section 3.4]{FeiZha} or    
\cite[Sections 3.2.2, 3.2.3]{FeiLuMal}.

We start with a test function 
\[
\phi \Grad \Del^{-1}[ \phi \beta (\vrd) ], \ \phi \in C^1_c (\Omega), 
\]
in the momentum balance \eqref{L3}:
\begin{equation} \label{L14}
\begin{split}
&\intTO{ \Big[ \vrd \vud \cdot \partial_t \left( \phi \Grad \Del^{-1}[ \phi \beta (\vrd) ] \right)  + \vrd \vud \otimes \vud : \Grad 
\left( \phi \Grad \Del^{-1}[ \phi \beta (\vrd) ] \right) \Big] }\\
&+ \intTO{ p_\delta (\vrd) \Div \left( \phi \Grad \Del^{-1}[ \phi \beta (\vrd) ] \right) }\\
&= \intTO{ \Big[ \mathbb{S}(\Grad \vud) : \Grad \left( \phi \Grad \Del^{-1}[ \phi \beta (\vrd) ] \right)  + \vrd \vc{g} \cdot \left( \phi \Grad \Del^{-1}[ \phi \beta (\vrd) ] \right)  \Big] }\\
&- \intTO{ \Lambda(\vrd) \vud \cdot \left( \phi \Grad \Del^{-1}[ \phi \beta (\vrd) ] \right) }.
\end{split}
\end{equation}

Next, using the renormalized equation \eqref{L15} we may identify the integral 
\begin{equation} \label{L16}
\begin{split}
&\intTO{ \vrd \vud \cdot \partial_t \left( \phi \Grad \Del^{-1}[ \phi \beta (\vrd) ] \right) } \\
&= - \intTO{ \vrd \vud \cdot  \left( \phi \Grad \Del^{-1}[ \phi \Div (\beta (\vrd) \vud ) ] ] \right) }\\
&+ \intTO{ \vrd \vud \cdot  \left( \phi \Grad \Del^{-1}[ \phi (\beta'(\vrd) \vrd - \beta (\vrd) )\Div \vud ) ] ] \right) }.
\end{split}
\end{equation}
Note that the boundary conditions are irrelevant here as $\phi$ is compactly supported. 

Now, we repeat the arguments of \cite{FeiZha} and consider in particular 
\[
\phi \Grad \Del^{-1} [ \phi \eta_\delta (\vrd) ]
\]
in \eqref{L14}, 
where 
\[
\eta_\delta (\vr) = \left\{ \begin{array}{l} \log(\Ov{\vr} - \vr) \ \mbox{if} \ 0 \leq \vr \leq \Ov{\vr} - \delta\\ \\
\log(\delta) \ \mbox{otherwise} \ \end{array} \right. .
\]
Given the available energy estimates \eqref{L9}, \eqref{L10}, and using the assumption in \eqref{R1} with $\beta\geq 3$ we deduce 
\begin{equation} \label{Pres1}
\int_{\mathcal{T}^1} \intO{ \phi^2 \eta_\delta (\vrd) p_\delta (\vrd) } \dt \leq c( \Grad \phi), 
\end{equation}
which yields the interior pressure bounds. Here it is hidden the role of the growth assumption on the pressure in \eqref{R1}, one can refer to \cite[Section 3.5]{FeiZha} for more details. 

To control the pressure up to the boundary, we repeat the process with 
the test function 
\[
\mathcal{B} [\phi] \ \mbox{for a suitable}\ \phi(x) \in W^{1,q}(\Omega; R^d), \phi(x) \to \infty \ \mbox{if}\ x \to \partial \Omega.
\]
If $q$ is large enough, we may use \eqref{Pres1} to get 
\[
\int_{\mathcal{T}^1} \intO{ \phi p_\delta (\vrd) } \dt, 
\]
which, together with \eqref{Pres1}, yields the desired equi--integrability of $\{ p_\delta (\vrd) \}_{\delta > 0}$, 
\begin{equation} \label{Pres2}
\int_{ p_\delta (\vrd) > \omega } p_\delta (\vrd) \dxdt \to 0 \ \mbox{for}\ \omega \to \infty.
\end{equation}
Thus 
$$ p_\delta (\vrd) \to \Ov{p(\vr)} \mbox{ weakly in } L^1((0, T)\times \Omega).$$
Finally, exactly as in \cite[Section 3.4]{FeiZha}, the choice 
\[
\phi \Grad \Del^{-1}[\phi \vrd], \ \ \phi  \in C^1_c (\Omega),
\]
gives rise to the so--called Lions identity, 
\begin{equation} \label{L17}
\Ov{p(\vr) \vr } - \Ov{p(\vr)} \vr = \left( \mu \left( 2 - \frac{2}{d} \right) + \eta \right) \left( \Ov{\vr \Div \vu} - \vr \Div \vu \right),
\end{equation}
or, more precisely, 
\[
\begin{split}
\lim_{\delta \to 0} &\intO{ \int_{\mathcal{T}^1} \phi^2 \Big[ p_\delta (\vrd) \vrd  - p_\delta (\vrd)  \vr \Big] } 
\dt\\ &= \lim_{\delta \to 0} \int_{ \mathcal{T}^1 } \intO{ \phi^2 \left( \mu \left( 2 - \frac{2}{d} \right) + \eta \right) \Big[ {\vrd \Div \vud} - \vr \Div \vu  \Big] } \dt. 
\end{split}
\]
As $p= p(\vr)$ is strictly increasing the left-hand side of \eqref{L17} is non-negative, thus relation \eqref{L17} implies convergence in measure (or a.a. convergence of a subsequence)
of $\{ \vrd \}_{\delta > 0}$, as soon as we show 
\begin{equation} \label{last}
\intTO{ \Ov{ \vr \Div \vu } - \vr \Div \vu } \leq 0,
\end{equation}
see \cite[Section 7]{CiFeJaPe1} for details.

\subsection{Proof of \eqref{last}   } 

Our ultimate goal is show \eqref{last}. In accordance with Remark \ref{renR},
we may consider  
\[
\varphi_\ep (t,x) = \min \left\{ 1 ; \frac{1}{\ep} {\rm dist}[(t,x); \Gamma_{\rm out}] \right\}
\]
as a test function in the renormalized equation of continuity \eqref{L15}, with $\beta(\vrd) = \vrd \log(\vrd)$. Performing the limit 
$\delta \to 0$ we obtain 
\begin{equation} \label{Pr1}
\intTO{ \left[ \Ov{\vr \log (\vr)} \partial_t \varphi_\ep + \Ov{ \vr \log(\vr)} \vu  \cdot \Grad \varphi_\ep - 
\Ov{ \vr \Div \vu } \varphi_\ep \right]  } = \int_{\Gamma_{\rm in}} 
\varphi_\ep \beta(\vr_b) \vuB \cdot \vc{n} \ \D \sigma_{x} \dt.
\end{equation} 
Applying the same treatment to the limit  equation, we get 
\begin{equation} \label{Pr2}
\intTO{ \left[ {\vr \log (\vr)} \partial_t \varphi_\ep + { \vr \log(\vr)} \vu  \cdot \Grad \varphi_\ep -
{ \vr \Div \vu } \varphi_\ep \right]  } = \int_{\Gamma_{\rm in}} 
\varphi_\ep \beta(\vr_B) \vuB \cdot \vc{n} \ \D \sigma_{x} \dt.
\end{equation}

Now, the task is to let $\ep \to 0$ in \eqref{Pr1}, \eqref{Pr2}. 
As $\partial \Omega$ is smooth (here we need at least $C^2$), there exists $\ep_0 > 0$ such that 
\begin{itemize}
\item
\[
\mathcal{U}_\ep (\partial \Omega) \equiv \left\{ x \in R^d \ \Big| \ {\rm dist}[x, \partial \Omega] < \ep \right\} 
= \cup_{x_b \in \partial \Omega}(x_b - \ep \vc{n}(x_b), x_b + \ep \vc{n}(x_b)),\ 0 < \ep < \ep_0; 
\]
\item 
the mapping 
\[
x_b \in \partial \Omega \mapsto x_b + \lambda \vc{n}(x_b) \ 
\ \mbox{is a diffeomorphism for any}\ \lambda \in [-\ep_0, \ep_0]
\]
\item for any $x \in \mathcal{U}_{\ep_0}$ there is a unique nearest point $x_b(x) \in \partial \Omega$, 
\[
\frac{x - x_b(x)}{|x - x_b(x)|} = - \vc{n}(x_b(x)), \ \mbox{if}\ x \in \mathcal{U}_{\ep_0} \setminus 
\partial \Omega.
\]

\end{itemize}

Next, we introduce the sets 
\[
\mathcal{O}_{\rm int} = \cup_{[t, x_b] \in \Gamma_{\rm int} } \left\{ [t, x_b + \lambda \vc{n}(x_b)] \Big| \ 
|\lambda| < \min \left\{ \ep _0 ; {\rm dist}[(t,x_b) ; \Gamma_{\rm out} ] \right\} \right\}, 
\]
and
\[
\mathcal{O}_{\rm out} = \cup_{[t, x_b] \in {\rm int}(\Gamma_{\rm out}) } \left\{ [t, x_b + \lambda \vc{n}(x_b)] \Big| \ 
|\lambda| < \min \left\{ \ep _0 ; {\rm dist}[(t,x_b) ; \Gamma_{\rm in} ] \right\} \right\} 
\]
The sets $\mathcal{O}_{\rm in}$, $\mathcal{O}_{\rm out}$ enjoy the following properties:
\begin{itemize}
\item 
$\mathcal{O}_{\rm in}$, $\mathcal{O}_{\rm out}$ are open in $\mathcal{T}^1 \times R^d$; 
\item 
\[
\Gamma_{\rm in} \subset \mathcal{O}_{\rm in},\ {\rm int}(\Gamma_{\rm out}) \subset \mathcal{O}_{\rm out}; 
\]
\item
\[
\begin{split}
{\rm dist} ([t,x] ; \Gamma_{\rm out}) &< \ep \ \mbox{for any}\ (t,x) \in \mathcal{O}_{\rm out},\ \mbox{whenever}\  
{\rm dist}[x, \partial \Omega] < \ep \\ &\Rightarrow \\ 
\varphi_\ep(t,x) &= \frac{1}{\ep} {\rm dist} ([t,x] ; \Gamma_{\rm out}),\ \partial_t \varphi_\ep (t,x) = 0, \ 
\Grad \varphi_\ep (t,x) = - \frac{1}{\ep} \vc{n}(x_b(x)),\ x_b(x) \in \Gamma_{{\rm out}};
\end{split}
\]
\item 
\[
\begin{split}
{\rm dist} ([t,x] ; \Gamma_{\rm out}) &> \ep \ \mbox{for any}\ (t,x) \in \mathcal{O}_{\rm in},\ \mbox{whenever}\ 
{\rm dist}[x, \partial \Omega] < \ep \\ &\Rightarrow \\ 
\varphi_\ep(t,x) &= 1,\ \partial_t \varphi_\ep (t,x) = 0, \ 
\Grad \varphi_\ep (t,x) = 0;
\end{split}
\]
\item 
\[
\frac{1}{\ep} \left| \left\{ (t,x) \in \mathcal{T}^1 \times \Omega \ \Big| \ 
{\rm dist}[x, \partial \Omega] \leq \ep, \ (t,x) \in (\mathcal{T}^1 \times \Omega) 
\setminus (\mathcal{O}_{\rm in} \cup \mathcal{O}_{\rm out}) \right\} \right| \to 0
\]
as $\ep \to 0$.

\end{itemize}

In view of the above observations, we may let $\ep \to 0$ in \eqref{Pr1} deducing
\begin{equation} \label{Pr3}
\begin{split}
\lim_{\ep \to 0}&
\intTO{ \left[ \Ov{\vr \log (\vr)} \partial_t \varphi_\ep + \Ov{ \vr \log(\vr)} \vu  \cdot \Grad \varphi_\ep - 
\Ov{ \vr \Div \vu } \varphi_\ep \right]  }\\&=
- \lim_{\ep \to 0} \frac{1}{\ep} \int_{\mathcal{O}_{{\rm out}} \cap \mathcal{U}_\ep (\partial \Omega)}\Ov{ \vr \log(\vr)} \vu (t,x)  \cdot \vc{n}(x_b(x)) \dxdt
\ \dxdt - \intTO{ \Ov{ \vr \Div \vu } }
 \\ &= \int_{\Gamma_{\rm in}}  \beta(\vr_b) \vuB \cdot \vc{n} \ \D \sigma_{x} \dt.
\end{split}
\end{equation}

Similarly, we get from \eqref{Pr2}  
\begin{equation} \label{Pr4}
\begin{split}
\lim_{\ep \to 0}&
\intTO{ \left[ {\vr \log (\vr)} \partial_t \varphi_\ep + { \vr \log(\vr)} \vu  \cdot \Grad \varphi_\ep - 
{ \vr \Div \vu } \varphi_\ep \right]  }\\&=
- \lim_{\ep \to 0} \frac{1}{\ep} \int_{\mathcal{O}_{{\rm out}} \cap \mathcal{U}_\ep (\partial \Omega)} { \vr \log(\vr)} \vu (t,x)  \cdot \vc{n}(x_b(x)) \dxdt
\ \dxdt - \intTO{ { \vr \Div \vu } }
 \\ &= \int_{\Gamma_{\rm in}}  \beta(\vr_b) \vuB \cdot \vc{n} \ \D \sigma_{x} \dt.
\end{split}
\end{equation}

Finally, we rewrite 
\begin{equation} \label{Pr5}
\begin{split}
\lim_{\ep \to 0} &\frac{1}{\ep} \int_{\mathcal{O}_{{\rm out}} \cap \mathcal{U}_\ep (\partial \Omega)}\Ov{ \vr \log(\vr)} \vu (t,x)  \cdot \vc{n}(x_b(x)) \dxdt\\ &= 
\lim_{\ep \to 0} \frac{1}{\ep} \int_{\mathcal{O}_{{\rm out}} \cap \mathcal{U}_\ep (\partial \Omega)}\Ov{ \vr \log(\vr)} (\vu (t,x) 
- \vuB (t,x))  \cdot \vc{n}(x_b(x)) \dxdt\\ 
&+ \lim_{\ep \to 0} \frac{1}{\ep} \int_{\mathcal{O}_{{\rm out}} \cap \mathcal{U}_\ep (\partial \Omega)}\Ov{ \vr \log(\vr)} (\vuB (t,x) 
- \vuB (t,x_b(x) ))  \cdot \vc{n}(x_b(x)) \dxdt\\ &+ 
\lim_{\ep \to 0} \frac{1}{\ep} \int_{\mathcal{O}_{{\rm out}} \cap \mathcal{U}_\ep (\partial \Omega)}\Ov{ \vr \log(\vr)} 
\vuB (t,x_b(x) )  \cdot \vc{n}(x_b(x)) \dxdt, 
\end{split}
\end{equation}
where 
\[
\begin{split}
\lim_{\ep \to 0} &\frac{1}{\ep} \int_{\mathcal{O}_{{\rm out}} \cap \mathcal{U}_\ep (\partial \Omega)}\Ov{ \vr \log(\vr)} (\vuB (t,x) 
- \vuB (t,x_b(x) ))  \cdot \vc{n}(x_b(x)) \dxdt\\ &+
\lim_{\ep \to 0} \frac{1}{\ep} \int_{\mathcal{O}_{{\rm out}} \cap \mathcal{U}_\ep (\partial \Omega)}\Ov{ \vr \log(\vr)} \frac{(\vu (t,x) 
- \vuB (t,x))}{{\rm dist}[x, \partial \Omega]} {\rm dist}[x, \partial \Omega]  \cdot \vc{n}(x_b(x)) \dxdt = 0,
\end{split}
\]
and, similarly, 
\[
\begin{split}
&\left| \lim_{\ep \to 0} \frac{1}{\ep} \int_{\mathcal{O}_{{\rm out}} \cap \mathcal{U}_\ep (\partial \Omega)}\Ov{ \vr \log(\vr)} (\vuB (t,x) 
- \vuB (t,x_b(x) ))  \cdot \vc{n}(x_b(x)) \dxdt \right| \\
& \leq \left\| \Grad \vuB \right\|_{L^\infty}  \lim_{\ep \to 0} \int_{\mathcal{O}_{{\rm out}} \cap \mathcal{U}_\ep (\partial \Omega)} |\Ov{ \vr \log(\vr)} | \dxdt = 0. 
\end{split}
\]

By the same token, we obtain 
\begin{equation} \label{Pr6}
\begin{split}
\lim_{\ep \to 0} &\frac{1}{\ep} \int_{\mathcal{O}_{{\rm out}} \cap \mathcal{U}_\ep (\partial \Omega)}{ \vr \log(\vr)} \vu (t,x)  \cdot \vc{n}(x_b(x)) \dxdt\\ = &\lim_{\ep \to 0} \frac{1}{\ep} \int_{\mathcal{O}_{{\rm out}} \cap \mathcal{U}_\ep (\partial \Omega)}{ \vr \log(\vr)} 
\vuB (t,x_b(x) )  \cdot \vc{n}(x_b(x)) \dxdt.
\end{split}
\end{equation}

Seeing that $\vr \mapsto \vr \log(\vr)$ is convex, and therefore $\Ov{\vr \log(\vr)} \geq \vr \log(\vr)$, and 
$\vuB \cdot \vc{n} \geq 0$ on $\Gamma_{\rm out}$, we combine \eqref{Pr4}--\eqref{Pr6} to obtain the desired conclusion 
\eqref{last}. We have proved Theorem \ref{wT1}.

\section{Concluding remarks}
\label{C} 

The result can be possibly extended to the full Navier--Stokes--Fourier system elaborating the arguments of \cite{FeMuNoPo}. Note, however, that the relevant existence theory in the case of the hard--sphere pressure for the evolutionary problem is to be developed. 
The stationary problem for $d=2$ and a tangential boundary velocity field has been treated in \cite{CiFeJaPe1}. 

The regularity of the boundary, and, in particular, its component $\Gamma_{\rm in}$ can be also relaxed, in the spirit of 
Chang et al. \cite{ChJiNo}. Similar extension in the case of Navier--Stokes--Fourier system is more delicate, see e.g. Poul \cite{Poul}. 

Last but not the least, the hard pressure enables to study the long--time behavior of the system, in particular the existence of bounded absorbing sets and attractors. 

\def\cprime{$'$} \def\ocirc#1{\ifmmode\setbox0=\hbox{$#1$}\dimen0=\ht0
  \advance\dimen0 by1pt\rlap{\hbox to\wd0{\hss\raise\dimen0
  \hbox{\hskip.2em$\scriptscriptstyle\circ$}\hss}}#1\else {\accent"17 #1}\fi}
\providecommand{\bysame}{\leavevmode\hbox to3em{\hrulefill}\thinspace}
\providecommand{\MR}{\relax\ifhmode\unskip\space\fi MR }
\providecommand{\MRhref}[2]{%
  \href{http://www.ams.org/mathscinet-getitem?mr=#1}{#2}
}
\providecommand{\href}[2]{#2}

\end{document}